\documentclass{article}

\usepackage[utf8]{inputenc}

\usepackage{amsmath}
\usepackage{amsfonts}
\usepackage{amssymb}

\usepackage{color}

\usepackage{siunitx}
\usepackage{graphicx}
\newcommand{\uv}{\mathbf{u}}
\newcommand{\nv}{\mathbf{n}}
\newcommand{\xv}{\mathbf{x}}
\newcommand{\fv}{\mathbf{f}}
\newcommand{\tv}{\mathbf{t}}
\newcommand{\ev}{\mathbf{e}}
\newcommand{\dv}{\mathbf{d}}
\newcommand{\vv}{\mathbf{v}}
\newcommand{\Ev}{\mathbf{E}}
\newcommand{\Dv}{\mathbf{D}}

\newcommand{\Piv}{\mathbf{P}^{i}}

\newcommand{\cten}{\mathbf{c}}

\newcommand{\hten}{\mathbf{h}}
\newcommand{\epsten}{\boldsymbol{\epsilon}}
\newcommand{\betaten}{\boldsymbol{\beta}}

\newcommand{\St}{\mathbf{S}}
\newcommand{\Tt}{\mathbf{T}}

\newcommand{\spaceV}{\mathbb{V}}
\newcommand{\spaceD}{\mathbb{D}}
\newcommand{\spaceP}{\mathbb{P}}
\newcommand{\spaceQ}{\mathbb{Q}}

\newcommand{\opdiv}{\operatorname{div}}

\newcommand{\citep}[1]{\cite{#1}}
\newcommand{\citet}[1]{\cite{#1}}

\begin{document}

\title{High-order mixed finite elements for an energy-based model of the polarization process in ferroelectric materials}

\author{Astrid S. Pechstein, Martin Meindlhumer and Alexander Humer}

\date{January 2020}
\maketitle

\begin{abstract}
	An energy-based model of the ferroelectric polarization process is presented in the current contribution. In an energy-based setting, dielectric displacement and strain (or displacement) are the primary independent unknowns. As an internal variable, the remanent polarization vector is chosen.  The model is then governed by two constitutive functions: the free energy function and the dissipation function. Choices for both functions are given. As the dissipation function for rate-independent response is non-differentiable, it is proposed to regularize the problem. Then, a variational equation can be posed, which is subsequently discretized using conforming finite elements for each quantity. We point out which kind of continuity is needed for each field (displacement, dielectric displacement and remanent polarization) is necessary to obtain a conforming method, and provide corresponding finite elements. The elements are chosen such that Gauss' law of zero charges is satisfied exactly. The discretized variational equations are solved for all unknowns at once in a single Newton iteration. We present numerical examples gained in the open source software package Netgen/NGSolve.
\end{abstract}

\section{Introduction}

%Piezoceramics are widely used as so-called smart materials for high-precision actuation and sensing. These abilities of piezoceramics are due to their ferroelectric behavior, i.e. their ability to switch their polarization state under sufficiently high external electric fields, and to maintain a remanent polarization as this poling field is removed.

Piezoceramics are widely used as so-called smart materials for high-precision actuation and sensing. 
Due to the piezoelectric effect, mechanical strains give rise to dielectric displacements, whereas applied voltages introduce stresses. To quantify the piezoelectric effect, the polarization state of the smart material has to be known. In most applications, the polarization state is assumed as constant unidirectional. This initial polarization is obtained by applying a high electric field to the ferroelectric specimen. Polarization aligned with the electric field builds up, and a remanent polarization is maintained as this poling field is removed.

However, the remanent polarization state can change under sufficiently high electrical or mechanical loadings. Also, in applications such as some macro-fibre composites (MFCs), the poling electric field is not unidirectional, which will lead to non-standard polarization patterns. Therefore, the modeling of the polarization process in ferroelectric materials and the subsequent numerical simulation are an active field of research.

More precisely, macroscopic modeling of the ferroelectric behavior has been pursued for a long time, as a macroscopic mathematical model is prerequisite to any simulation using finite elements or other concepts. The model is required to catch the nonlinear material behavior especially under high electric fields and also depolarization under high pressures. Two inherently different approaches have been pursued successfully so far: microscopic-based and phenomenological models.

Microscopic-based models work on the level of microscopic cells or crystals, from which a macroscopic behaviour is derived by averaging over a large number of unit cells. This approach has been discussed e.g.\ by Huber et al.\ \cite{HuberFleckLandisMcMeeking:1998}. Finite element methods based on microscopic models date back to Chen and Lynch \cite{ChenLynch:1998}. With increasing computing power, methods using one finite element per grain have become feasible, we cite non-exhaustively \cite{ArockiarajanSivakumarSansour:2009,LangeRicoeur:2015,JayendiranGanapathiZineb:2016,LobanovSemenov:2019}.

Another approach, which is adopted in this work, relies on a phenomenological description of the macroscopic material behavior. Internal variables then resemble macroscopic quantities, such as the remanent polarization or remanent polarization strain. We differ between thermodynamically consistent and non-consistent models.

Models based on a Preisach operator are not thermodynamically consistent. Parameters are chosen to emulate hysteretic curves. This approach has been chosen e.g.\ by \cite{HwangEtal:1995,HughesWen:1995,HegewaldEtal:2008}. Non-linear finite elements for the latter model have been proposed by \cite{KKaltenbacherHegewaldLerch:2010}.

A thermodynamically consistent model was first presented in the series of papers by the group around Maugin \citep{BassiounyGhalebMaugin:1988a,BassiounyGhalebMaugin:1988b,BassiounyMaugin:1989a,BassiounyMaugin:1989b}. Kamlah \cite{Kamlah:2001} developed a theory for multi-axial loading based on several switching criteria defining the polarization and strain hystereses. A finite element implementation was proposed by Kamlah and Böhle  \cite{KamlahBoehle:2001}. More recently, Elhadrouz et al.\ \cite{ElhadrouzZinebPatoor:2005} proposed a model including two different switching criteria for ferroelectric and ferroelastic switching.

Another family of thermodynamically consistent models has been developed, where ferroelectric and ferroelastic switching are combined in only one switching function. McMeeking and Landis \cite{McMeekingLandis:2002} as well as Schröder and Romanowski \cite{SchroederRomanowski:2005} proposed formulations where the remanent polarization strain depends one-to-one on the remanent polarization. Nevertheless, these models are capable of mechanic depolarization, see the numerical examples from \cite{McMeekingLandis:2002}. Linnemann et al.\ \cite{LinnemannKlinkelWagner:2009} added magnetostriction in their work. Another approach with an independent polarization strain has been proposed by Landis \cite{Landis:2002}. Also Klinkel \cite{Klinkel:2006} uses two independent remanent quantities remanent electric field and remanent strain in his model. Laxman et al.\ \cite{LaxmanManiprakashArockiarajan:2018} additionally modeled creep behavior  in their approach.

Miehe et al.\ \cite{MieheRosatoKiefer:2011} presented a variational framework for electro-magneto-mechanics. They introduced a concept of minimizing the work increment. They differ between approaches based on the free energy and electric enthalpy based approaches. Depending on this base, the free unknowns are either strain and dielectric displacement, or strain and electric field. While all of the models mentioned above are electric enthalpy based, only few energy-based models have been used so far. Sands and Guz \cite{SandsGuz:2013} describe a one-dimensional model based on the free energy. Semenov et al.\ \cite{SemenovLiskowskyBalke:2010} provide a multi-directional energy based model using a vector potential for the dielectric displacement. They propose to discretize the vector potential in the finite element method, and provide consistent tangent moduli for the implementation. They observed very fast convergence of the Newton-Raphson iteration. Stark et al.\ \cite{StarkNeumeisterBalke:2016I,StarkNeumeisterBalke:2016II} combine a phenomenological and micro-mechanically motivated model in their work.

The authors of the present work have analyzed such an energy-based formulation in the framework of variational inequalities in \cite{PechsteinMeindlhumerHumer:2019arxiv}. This concept has been used in the mathematical analysis of other nonlinear problems in mechanics, such as elasto-plasticity by Han and Reddy \cite{HanReddy:1999} or contact mechanics by Sofonea and Matei \cite{SofoneaMatei:2011}. In the present work, a finite element method for the energy-based formulation is proposed. The utilized finite elements for displacement, dielectric displacement and remanent polarization are described in detail. These elements are chosen such that Gauss' law of zero charges inside the domain is satisfied exactly, what is usually not the case in electric enthalpy based methods. Numerical results for the problem of non-proportional loading and the homogenization of the polarization of a macro-fibre composite are presented.

The paper is organized as follows: First, the energy-based model in analogy to \cite{MieheRosatoKiefer:2011} is presented. Algorithmic representations for energy and dissipation are provided subsequently. A mathematical model of the time-discrete polarization problem is provided, which is non-differentiable due to the choice of dissipation function. Its regularization is addressed, and also the stability of the resulting variational equation. The next section is devoted to the finite element discretization, the choice of finite elements and the implementation in the open-source software package Netgen/NGSolve. The last section comprises numerical results.

\section{Modeling dissipative material behavior}

In the following, we assume to be in the standard setting of small-deformation mechanics in three space dimensions. The body of interest shall be denoted by $\Omega$, which is in total or in part made of ferroelectric material. We use $\uv: \Omega \times [0,T] \to \mathbb{R}^3$ to denote the displacement field, and use $\St = \frac12 (\nabla \uv + (\nabla \uv)^T)$ for the linearized strain tensor. 

In an electrostatic regime on a simply connected domain, the electric field $\Ev$ can be shown to be a gradient field due to Faraday's law. We use standard nomenclature from the literature, introducing the electric potential $\varphi$ with $\Ev =  - \nabla \varphi$. 
The dielectric displacement vector $\Dv$ satisfies Gauss' law of zero charge in the material,
\begin{align} \label{eq:Gauss}
\opdiv \Dv &= 0 &&\text{in } \Omega.
\end{align}
Last, we use the remanent polarization $\Piv$ as an internal variable.

The reversible material response is characterized by the internal energy $\Psi$ or its density $\psi$, which shall depend on strain, dielectric displacement as well as remanent polarization vector,
\begin{align}
\psi(\xv, t) &= \hat \psi(\St(\xv,t), \Dv(\xv,t), \Piv(\xv,t)),\\
\Psi(t) &= \int_\Omega \psi(\xv,t)\, d\xv.
\end{align}

Total stress $\Tt$ and electric field $\Ev$ are conjugate to strain $\St$ and dielectric displacement $\Dv$ in the following sense,
\begin{align}
\Tt &= \frac{\partial \hat\psi}{\partial \St}, & \Ev &= \frac{\partial \hat \psi}{\partial \Dv}.
\label{eq:defTE}
\end{align}
Conjugate to the remanent polarization is the (negative) electric driving force $\hat \Ev$,
\begin{align}
\hat \Ev &= -\frac{\partial \hat \psi}{\partial \Piv}.
\label{eq:defhatE}
\end{align}

We denote by $P^{ext} = \frac{d}{dt} W^{ext}$ the power added from external sources such as applied forces or electric fields. As in mixed methods the external work often differs from standard approaches, we provide an exemplary expression below. For a given volume load $\fv$ in $\Omega$,   a surface load $\tv$ on $\Gamma_{t}$  and an applied electric potential $V_0$ on $\Gamma_{el}$, the power from external sources is given by
\begin{align} \label{eq:externalwork}
P^{ext} &= \int_\Omega \fv\cdot \dot\uv\, d\xv + \int_{\Gamma_{t}} \tv \cdot \dot\uv\, ds - \int_{\Gamma_{el}} V_0\, \dot\Dv\cdot \nv\, ds.
\end{align}

A reversible process is fully characterized by internal energy $\Psi$ and external work $W^{ext}$, and there are no internal variables such as the remanent polarization present. In this contribution however, we are concerned with irreversible, dissipative processes. 
As stated by Miehe et al.\ \cite{MieheRosatoKiefer:2011}, the dissipative behavior can be described by the dissipation function $\Phi$, which depends only on the rate of the internal variable. In our case, the dissipation function is thus a function of the remanent polarization rate,
\begin{align}
\phi(\xv,t) &= \hat \phi(\dot \Piv(\xv,t)),\\
\Phi(t) &= \int_\Omega \phi(\xv,t)\, d\xv.
\end{align}
Depending on the dissipation function, rate-dependent or rate-independent material behavior is modeled. In the next section, we motivate our special choice of dissipation function, which leads to rate-independent behavior.

Using the notions introduced above, we cite the following incremental minimization problem for any time interval $(t_0, t_1)$,
\begin{align}
\begin{split}
\Psi(t_1) - \Psi(t_0)& + \int_{t_0}^{t_1}\!\! \Phi(t)\, dt - \int_{t_0}^{t_1}\!\!\! P^{ext}(t)\, dt \to \min_{\uv, \Dv}\\
\text{subject to } & \left\{\begin{array}{l}\opdiv \Dv = 0 \text{ and}\\
\text{boundary conditions on } \uv, \Dv\cdot\nv.\end{array} \right.
\end{split} \label{eq:minimization_cont}
\end{align}
For a more detailed introduction, we refer to the extensive derivations originally provided by \cite{MieheRosatoKiefer:2011}.

\subsection{Choice of internal energy}

Different choices of the internal energy density have been proposed in the literature. 
In many contributions, the internal energy is decomposed into two additive terms: the first of these terms coincides with the energy in linear problems for fixed remanent polarization. It is quadratic in the reversible parts of dielectric displacement and strain, but material tensors $\cten^{D}, \hten$ and $\betaten^{S}$ may depend on the remanent polarization. The second term depends then only on the remanent polarization, and contains the electric hystersis shape parameters $h_0$ and $m$ and the saturation polarization $P_0$.
\begin{align}
\begin{split}
\hat \psi =& \ \hat \psi^{r} +\hat  \psi^{i},\\
\hat \psi^r =&\ \frac12 (\St - \St^{i}(\Piv)):\cten^{D}:(\St - \St^{i}(\Piv))\\
& -(\St - \St^{i}(\Piv)):\hten\cdot(\Dv - \Piv) \\
& +\frac12 (\Dv - \Piv)\cdot\betaten^{S}\cdot(\Dv-\Piv),\\
\hat \psi^{i} =&\ \frac{h_0\, P_0^2}{(m-1)(m-2)} \left( 1-\frac{|\Piv|}{P_0} \right)^{2-m}\\
& - \frac{h_0\, P_0^2}{m-1}\frac{|\Piv|}{P_0}.
\end{split} \label{eq:defpsi}
\end{align}
We use a one-to-one relationship between remanent polarization and polarization strain as suggested by \cite{McMeekingLandis:2002},
\begin{align}
S^{i}_{kl}(\Piv) &= \frac{S_0}{2 P_0^2} (3 P^{i}_k P^{i}_l - \delta_{kl} P^{i}_j P^{i}_j).
\end{align}
Note that this choice is volume preserving.

\subsection{Choice of dissipation function}
Apart from the free energy density $\hat \psi$, the dissipation function $\hat \phi$ has to be specified to obtain a valid model through the minimization problem \eqref{eq:minimization_cont}. The dissipation function is closely related to the dissipation $D$. The dissipation describes the amount of energy dissipated into heat and is always non-negative, as stated by the Clausius-Duhem inequality
\begin{align}
D = \dot \Piv \cdot \hat \Ev \geq 0.
\end{align}
The choice of dissipation function $\hat \phi$ is motivated by the principle of maximum dissipation. It is closely related to the switching condition, which states in its simplest form that the remanent polarization rate is non-trivial only if the driving electric field $\hat \Ev$ exceeds the coercive field $E_0$,
\begin{align}
|\hat \Ev| - E_0 & \leq 0, & |\dot \Piv| &\geq 0, & (|\hat \Ev| - E_0)|\dot \Piv| = 0.
\end{align}
The dissipation function is chosen such that dissipation is maximized over all possible values for $\hat \Ev$, 
\begin{align}
\hat \phi(\dot \Piv) = \sup_{\hat \Ev, |\hat \Ev|\leq E_0} \dot \Piv\cdot \hat \Ev.
\end{align}
This supremum is of course taken for $\hat \Ev = E_0 \ev_P$, with $\ev_P = \Piv/|\Piv|$ the unit vector in direction of remanent polarization, and for $\hat \Ev = \mathbf{0}$ in case $|\dot \Piv| = 0$. Thus, the dissipation function can be simplified to read
\begin{align}
\hat \phi(\dot \Piv) = E_0 |\dot \Piv|.
\end{align}
Note that this dissipation function is positively homogeneous (compare \eqref{eq:poshom}). Such a dissipation function leads to rate-independent behavior, see \citet{MieheRosatoKiefer:2011}.
\section{Mathematical modeling}

To get a well-defined problem in a mathematical sense, we introduce Hilbert spaces for the various quantities. Then, the time-discrete minimization problem \eqref{eq:minimization_cont} can be solved minimizing in the respective closed spaces. We use the standard $H^1$ Sobolev space of weakly differentiable functions with derivative in $L^2$ for the displacements. (Homogeneous) displacement boundary conditions on the boundary part $\Gamma_{fix}$ are included, such that we are searching for
\begin{align}
\uv \in \spaceV &:= \{\vv \in H^1(\Omega): \vv = 0 \text{ on } \Gamma_{fix}\}.
\end{align}
The choice of the dielectric displacement space is less standard. We resort to the $H(\opdiv)$ space of $L^2$ vector fields with weak divergence. Gauss' law \eqref{eq:Gauss} is then included explicitely. Additionally, insulation boundary conditions of vanishing surface charges $\Dv \cdot \nv = 0$ are posed explicitely. This means that these surface charge boundary conditions are essential boundary conditions, while we have already seen in the power statement \eqref{eq:externalwork} that voltage boundary conditions on electrodes are natural. Summing up, we search for 
\begin{align}
\begin{split}
\Dv \in \spaceD_0 &:= \{\dv \in \spaceD: \opdiv \dv = 0\},\\
\text{with } \spaceD & = \{\dv \in [L^2]^3: \opdiv \dv \in L^2,\ \dv\!\cdot\!\nv = 0 \text{ on } \Gamma_{ins}\}\\
& \subset H(\opdiv).
\end{split}
\end{align}

The remanent polarization does not sport any smoothness, nor does it allow for boundary conditions. Therefore we choose the vector-valued Lebesgue space $[L^2]^3$,
\begin{align}
\Piv \in \spaceP &:= [L^2]^3.
\end{align}

We will now formulate the incremental variational problem for some small but finite time increment $\Delta t$, starting from time $t_0$ to $t_1 = t_0+\Delta t$. To this end, we assume $\uv_0 = \uv(t_0)$, $\Dv_0 = \Dv(t_0)$ and $\Piv_0 = \Piv(0)$ known. We are interested in finding  the increments $\Delta \uv, \Delta \Dv$ and $\Delta \Piv$ such that $\uv(t_1) = \uv_0 + \Delta \uv$, $\Dv(t_1) = \Dv_0 + \Delta \Dv$ and $\Piv(t_1) = \Piv_0 + \Delta \Piv$. Miehe et al.\ \cite{MieheRosatoKiefer:2011} stated that these increments minimize the potential increment,
\begin{align}
\begin{split}
&(\Delta \uv, \Delta \Dv, \Delta \Piv) =\\
&= \arg\min \Big( \Psi(t_1) - \Psi(t_0) + \Delta t \Phi(\Delta \Piv/\Delta t) - W^{ext}_{t_0\to t_1}\Big)
\\
& \quad \text{subject to } \Delta \uv \in \spaceV, \Delta \Dv \in \spaceD_0, \Delta \Piv \in \spaceP.
\end{split}
\label{eq:minimization_deltat_ratedep}
\end{align}
In the special case of a rate-independent response, where the dissipation function is positively homogeneous,
\begin{align} \label{eq:poshom}
\Phi(\alpha \dot \Piv) &= \alpha \Phi(\dot \Piv) & \text{for } \alpha > 0,
\end{align}
the time step size cancels out in \eqref{eq:minimization_deltat_ratedep}. The parameter $t$ acts as pseudo-time then, the minimization problem reads
\begin{align}
\begin{split}
&(\Delta \uv, \Delta \Dv, \Delta \Piv) =\\
&= \arg\min \Big( \Psi(t_1) - \Psi(t_0) + \Phi(\Delta \Piv) - W^{ext}_{t_0\to t_1}\Big)
\\
&\quad \text{subject to }   \Delta \uv \in \spaceV, \Delta \Dv \in \spaceD_0, \Delta \Piv \in \spaceP.
\end{split}
\label{eq:minimization_deltat}
\end{align}

\subsection{Variational inequality}
As the dissipation function is not differentiable at $\dot \Piv = \mathbf{0}$, this minimization problem cannot be transformed into a variational equation. However, it can be treated in the framework of variational inequalities. Inspired by the mathematical analysis of elasto-plasticity (see e.g. \cite{HanReddy:1999}) or contact problems (see e.g. \cite{SofoneaMatei:2011}), a thorough analysis of the arising time-dependent variational inequality has been carried out in \cite{PechsteinMeindlhumerHumer:2019arxiv}. In this contribution, we state the variational inequality which is equivalent to the optimization problem \eqref{eq:minimization_deltat}. We use the definitions for total stress $\Tt$ and electric field $\Ev$ \eqref{eq:defTE} and for the electric driving force $\hat \Ev$ \eqref{eq:defhatE} evaluated at $t_1$,
\begin{align}
\Tt_{t_1} &= \frac{\partial \hat\psi}{\partial \St}(\St(\uv_0+\Delta \uv),\Dv_0 + \Delta \Dv, \Piv_0 + \Delta \Piv),\\
\Ev_{t_1} &= \frac{\partial \hat\psi}{\partial \Dv}(\St(\uv_0+\Delta \uv),\Dv_0 + \Delta \Dv, \Piv_0 + \Delta \Piv),\\
\hat \Ev_{t_1} &= -\frac{\partial \hat\psi}{\partial \Piv}(\St(\uv_0+\Delta \uv),\Dv_0 + \Delta \Dv, \Piv_0 + \Delta \Piv).
\end{align}
Then, $\Delta \uv \in \spaceV$, $\Delta \Dv \in \spaceD_0$ and $\Delta \Piv \in \spaceP$ satisfy the following variational inequality
\begin{align}
\begin{split}
\int_\Omega \big( \Tt_{t_1}:(\St(\tilde \uv)-\St(\Delta\uv)) + \Ev_{t_1} \cdot (\tilde \Dv - \Delta \Dv) &\\
- \hat \Ev_{t_1} \cdot ( \tilde \Piv - \Delta \Piv)
+ \hat \phi(\tilde \Piv) - \hat \phi(\Delta \Piv) \big) d\xv &\\
\geq  \tilde W^{ext}_{t_0\to t_1} - W^{ext}_{t_0\to t_1} &\\
\text{for all } \tilde \uv \in \spaceV, \tilde \Dv \in \spaceD_0, \tilde \Piv \in \spaceP.&
\end{split} \label{eq:varin}
\end{align}

\subsection{Regularization}
In the framework of elasto-plasticity, Han and Reddy \cite{HanReddy:1999} propose to regularize the dissipation function. The dissipation function is changed slightly, such that it becomes differentiable everywhere. When using the regularized dissipation function, the variational inequality \eqref{eq:varin} turns into a nonlinear variational equation, which can be solved by standard nonlinear solution algorithms.

We propose to use the regularized dissipation function
\begin{align}
\hat \phi_\varepsilon &= \left\{ \begin{array}{ll}
	E_0\big(|\dot\Piv| - \frac{\varepsilon}{2}\big) & |\dot \Piv| \geq \varepsilon\\
	\frac1{2\varepsilon} E_0 |\dot \Piv|^2 & |\dot \Piv| < \varepsilon.
	\end{array} \right.
	\label{eq:defphieps}
\end{align}
This choice has been analyzed to work well for elasto-plasticity in \cite{HanReddy:1999}. Indeed, according to this reference, we expect the error due to regularization to be of order $\sqrt{\varepsilon}$.

For completeness we provide the variational equation to be solved in each time step. We denote the (continuous) derivative of the regularized dissipation function $\hat \phi_\varepsilon$ with respect to the remanent polarization rate as $\hat \phi_\varepsilon'$.
The updates $\Delta \uv \in \spaceV$, $\Delta \Dv \in \spaceD_0$ and $\Delta \Piv \in \spaceP$ satisfy 
\begin{align}
\begin{split}
\int_\Omega \big( \Tt_{t_1}:\,\delta \St + \Ev_{t_1} \cdot \delta \Dv - \hat \Ev_{t_1} \cdot \delta \Piv&\\
+ \hat \phi'_\varepsilon(\Delta \Piv)\cdot \delta \Piv  \big) &d\xv
=  \delta W^{ext}_{t_0\to t_1}  \\
\text{for all } \tilde \uv \in \spaceV, \tilde \Dv \in \spaceD_0, \tilde \Piv \in \spaceP.&
\end{split} \label{eq:vareq}
\end{align}

Another issue for regularization is the free energy function $\psi$, or its additive part $\psi^{i}$. The formula proposed in \eqref{eq:defpsi} is unbounded as the saturation polarization is approached,
\begin{align}
\hat \psi^{i}(\Piv) \to \infty &\text{ and } (\hat \psi^{i})'(\Piv)\cdot \delta\Piv \to \infty \\
&\text{ as } |\Piv| \to P_0.
\end{align}
This unbounded potential acts as a penalty term as the polarization approaches saturation, where also $ (\hat \psi^{i})'(\Piv)$  diverges to infinity.
As analyzed in \cite{PechsteinMeindlhumerHumer:2019arxiv}, this fact does not harm solvability of the update equation \eqref{eq:vareq} itself, but it may harm convergence. Therefore, if convergence problems are observed, we propose to modify the energy density $\hat \psi^{i}$ in a way that $ (\hat \psi^{i}_\varepsilon)'(\Piv)$  becomes large but stays bounded. We use the following regularized expression $ (\hat \psi^{i}_\varepsilon)'$,
\begin{align}
(\hat \psi^{i}_\varepsilon)'(\Piv) &= \left\{ \begin{array}{l}
\left( \frac{h_0 P_0^m}{m-1}(P_0-|\Piv|)^{1-m} -1\right) \Piv/|\Piv| \\
\text{  if } |\Piv| \leq P_0-\varepsilon,\\
\left(\frac{h_0 P_0^m}{m-1}(\varepsilon)^{1-m} -1\right) \Piv/(P_0-\varepsilon) \\
\text{ if } |\Piv| > P_0-\varepsilon.
\end{array} \right. 
\end{align}
\subsection{Gauss' law}

Gauss' law of zero charges \eqref{eq:Gauss} implies that the dielectric displacement field is divergence free inside the domain $\Omega$. So far, this condition has been included in the space of admissible updates $\spaceD_0$ without further ceremony. As Gauss' law is a linear and continuous restriction, the space $\spaceD_0$ is a linear and continuous subspace of $\spaceD$ and $H(\opdiv)$. However, as this condition is non-local, it will be difficult to implement a finite element space consisting of divergence-free finite element functions in practice. In the following paragraph, we discuss how to reformulate the condition in an equivalent way, such that finit element implementation will be possible. We stress that, both in theory in this section, as well as in the finite element implementation, the dielectric displacement field will be exactly divergence free in the strong sense.

We propose to add a Lagrangian multiplier ensuring the divergence-free condition. As it turns out that this multiplier resembles the electric potential, we denominate it as $\varphi$. In the minimization problem, we propose to use an augmented energy $\Psi_0$, which assumes infinity if Gauss' law is violated. This augmented energy reads
\begin{align}
\Psi_0 &= \Psi(\St, \Dv, \Piv) - \inf_{\varphi \in \spaceQ} \int_\Omega \opdiv \Dv\, \varphi\, d\xv.
\end{align}
Note that we use $\varphi \in \spaceQ := L^2(\Omega)$, which means that discontinuous electric potentials are admissible testing functions. Also, the electric potential does not satisfy any boundary conditions a priori. Of course, as the approach is consistent, the potential will assume boundary values matching the prescribed external work \eqref{eq:externalwork}. Also, it is not necessary to implement an equi-potential condition on electrodes specifically.

The regularized variational equation including Gauss' law is stated below. 
The updates $\Delta \uv \in \spaceV$, $\Delta \Dv \in \spaceD$ and $\Delta \Piv \in \spaceP$ as well as the (total) electric potential $\varphi \in \spaceQ$ satisfy 
\begin{align}
\begin{split}
\int_\Omega \big( \Tt_{t_1}:\,\delta \St + \Ev_{t_1} \cdot \delta \Dv - \hat \Ev_{t_1} \cdot \delta \Piv&\\
-\opdiv(\delta \Dv)\, \varphi - \opdiv (\Delta \Dv) \, \delta \varphi&\\
+ \hat \phi'_\varepsilon(\Delta \Piv)\cdot \delta \Piv & \big) d\xv
=  \delta W^{ext}_{t_0\to t_1}  \\
\text{for all } \delta \uv \in \spaceV, \delta \Dv \in \spaceD, \delta \Piv \in \spaceP, \delta \varphi \in \spaceQ.&
\end{split} \label{eq:vareq_lag}
\end{align}

The augmented variational equation \eqref{eq:vareq_lag} is equivalent to \eqref{eq:vareq}. The dielectric displacement update $\Delta \Dv$ is exactly divergence free and thus $\Delta \Dv \in \spaceD_0$.

\subsection{Solvability and stability}

The issue of existence of solutions, also to the time-dependent variational inequality \eqref{eq:varin}, has been treated by the authors in \citet{PechsteinMeindlhumerHumer:2019arxiv}. We recall that a unique solution to the update problem \eqref{eq:varin} exists, given the following essential assumption of stong monotonicity holds: for all $\uv, \tilde \uv \in \spaceV$, $\Dv, \tilde \Dv \in \spaceD_0$ and $\Piv, \tilde \Piv \in \spaceP$ as well as corresponding forces $\Tt, \tilde \Tt$, $\Ev, \tilde \Ev$ and $\hat \Ev, \tilde{\hat \Ev}$ there holds
\begin{align}
\int_\Omega \big( &(\Tt - \tilde \Tt):(\St - \tilde \St) + (\Ev - \tilde \Ev)\cdot(\Dv - \tilde \Dv) \nonumber\\
& - (\hat \Ev - \tilde {\hat \Ev})\cdot(\Piv - \tilde \Piv) \big)\, d\xv \\
&\geq C \big( \|\uv - \tilde \uv\|_\spaceV^2 + \|\Dv - \tilde \Dv\|_\spaceD^2 + \|\Piv-\tilde \Piv\|_\spaceP^2\big).\nonumber
\end{align}
This is the case if the reversible part of energy $\Psi^{r}$ is strongly monotone with respect to the reversible parts $\St - \St^{i}(\Piv)$ and $\Dv - \Piv$, and the additive part $\Psi^{i}$ is strongly monotone with respect to $\Piv$. This is the case, as long as the linear moduli lead to a positive definite material tensor. This issue is addressed by \cite{StarkNeumeisterBalke:2016a}. Conditions ensuring positive definiteness of the additive energy part $\hat \psi^{i}$ can be found in \cite{BotteroIdiart:2018}.

\section{Finite elements}

In this section we propose a set of consistent finite elements to discretize the variational equation \eqref{eq:vareq_lag} above. The elements shall be of arbitrary order, by $k \geq 1$ we denote the lowest overall polynomial order. In the following, we assume that $\mathcal{T}$ is a simplicial mesh of the volume $\Omega$. The elements can be defined also on prismatic meshes, as is shown in our numerical results.
 Note that, in the context of our mixed elements, ``consistent'' means different notions of continuity for the different fields.
 \subsection{Displacement elements}
  The displacement element has to be consistent for $H^1$, such that $\uv \in \spaceV$. This condition is satisfied when using standard nodal or hierarchical continuous elements. We propose to use hierarchical elements of arbitrary order. For simplicial elements, the finite element space $\spaceV_h$ is then given by
 \begin{align}
 \spaceV_h &:= \{\vv \in \mathcal{C}(\Omega): \vv|_T \in [P^{k+1}]^3\text{ for } T \in \mathcal{T}\} \subset \spaceV.
 \end{align}
 For tetrahedral elements, this results in $(k+2)(k+3)(k+6)/2$ degrees of freedom, which means $30$ degrees of freedom for $k=1$ (second order displacment elements) and $60$ degrees of freedom for $k=2$ (third order displacment element). A prismatic element sports $3(k+1)^2(k+2)/2$ degrees of freedom. For $k=1$ this results in $54$ degrees of freedom (second order displacment elements), for $k=2$ we obtain $120$ degrees of freedom (third order displacment elements).
 
 \subsection{Dielectric displacement and electric potential}
 The dielectric displacement elements are less standard. To be consistent, we need that $\Dv \in \spaceD$, i.e. that $\Dv$ allows for a divergence in weak sense. This means that the elements need to be normal continuous. In other words, different moments of the normal component $\Dv \cdot \nv$ on element interfaces are degrees of freedom instead of nodal values. Different families of such normal-continuous elements have been introduced in the literature. For a fundamental introduction and overview, we refer the interested reader to the monograph \cite{BoffiBrezziFortin:2013}. As a first proposition, we propose to use elements from the BDM family as introduced by \cite{BrezziDouglasDuranFortin:1987} (BDM after Brezzi, Douglas and Marini, who first proposed corresponding elements in two dimensions, see \cite{BrezziDouglasMarini:1986}). These elements span the space

 \begin{align}
 \begin{split}
 \spaceD_h:= \{ &\Dv \in [L^2(\Omega)]^3: \Dv|_T \in [P^k]^3\text{ for } T \in \mathcal{T},\\
 & \Dv \cdot \nv \text{ cont.}, \Dv \cdot \nv = 0 \text{ on }\Gamma_{ins}\} \subset \spaceD.
 \end{split}
 \end{align}
 
 However, as the dielectric displacement is divergence free, the number of unknowns can be further reduced. As already mentioned, the divergence-free condition is non-local, therefore it is not possible to restrict the space $\spaceD_h$ such that it becomes a subspace of $\spaceD_0$ in a simple, local manner. However, in context of incompressible Navier-Stokes equation, \cite{LehrenfeldSchoeberl:2016} have shown that it is possible to do so for the higher-order shape functions. We propose to use this reduced set,
 \begin{align}
\begin{split}
\spaceD_{h0}:= \{ &\Dv \in [L^2(\Omega)]^3: \Dv|_T \in [P^k]^3\text{ for } T \in \mathcal{T},\\
& \opdiv \Dv|_T \in P^0 \text{ for } T \in \mathcal{T},\\
& \Dv \cdot \nv \text{ cont.}, \Dv \cdot \nv = 0 \text{ on }\Gamma_{ins}\} \subset \spaceD.
\end{split}
\end{align}
 
 We see that this way not only the number of unknowns for $\Dv$ drops, but also less Lagrangian multipliers $\varphi$ are needed. As $\opdiv \Dv$ is constant per element, it is enough to use piecewise constant $\varphi$,
 \begin{align}
 \spaceQ_h &:= \{ \varphi \in L^2(\Omega): \varphi|_T \in P^0\text{ for } T \in \mathcal{T}\}
 \end{align}
 Then, the condition 
 \begin{align}
 \int_\Omega \opdiv \Dv \, \delta \varphi \, &= 0 && \text{for all } \delta \varphi \in \spaceQ_h,
 \end{align}
 implies that $\opdiv \Dv = 0$ in strong sense for $\Dv \in \spaceD_{h0}$.
 
 The finite element space $\spaceD_{h0}$ is intrinsically different from standard nodal finite element spaces well-known in mechanics. While in nodal finite element spaces such as $\spaceV_h$ degrees of freedom are function values evaluated at element nodes, a more abstract concept of degrees of freedom as introduced by \cite{Ciarlet:1978} is needed in the context of $H(\opdiv)$. In our case, degrees of freedom are associated to element interfaces (to $\Dv\cdot \nv$), and additional ``element bubbles'' that are associated to the element interior (similar to element-interior nodes). There are no degrees of freedom associated to nodes or element edges. Different from standard nodal elements as used for the displacements, the inter-element coupling acts only through degrees of freedom associated to element faces. Therefore, the coupling is less strong than for standard elements, and even for the same number of overall degrees of freedom, the system matrix is sparser and thereby easier to invert by a direct solver. As the system matrix is indefinite, it is important to choose a direct solver that can handle such systems.
 
 We count the total number of degrees of freedom for the electric quantities per element. In all cases, one degree of freedom is used for the electric potential, while the number of dielectric displacement degrees of freedom varies depending on $k$ and the element shape. For a tetrahedral element of order $k$, there are $2(k+1)(k+2)$ coupling and $k(k-1)(2k+5)/6$ interior degrees of freedom. For $k=1$ this means $12$ coupling degrees of freedom only, for $k=2$ this adds up to $24$ coupling and $3$ interior degrees of freedom. For a prismatic element, we have $(k+1)(k+2)+3(k+1)^2$ coupling and $(k+2)(k+1)(2k+1)/2+1$ interior degrees of freedom. For a first order element with $k=1$ this adds up to $18$ coupling and $10$ inner degrees of freedom. For a second order element with $k=2$ we get $39$ coupling and $31$ inner degrees of freedom.
 
 \subsection{Remanent polarization}
 There are no conditions on the remanent polarization other than $\Piv \in \spaceP = [L^2]^3$. Neither in derivation nor in stability analysis any assumptions on the continuity of $\Piv$ were made. Therefore, we propose the remanent polarization to be of the same order as the dielectric displacement and discontinuous,
  \begin{align}
 \spaceP_h &:= \{ \Piv \in [L^2(\Omega)]^3: \Piv|_T \in [P^k]^3\text{ for } T \in \mathcal{T}\}
 \end{align}
 
\subsection{Implementation}

All finite elements described in this section are implemented in the open source software package Netgen/NGSolve available from {\tt ngsolve.org}. Netgen/NGSolve sports a python interface, where one can symbolically enter variational equations or energies. The necessary derivations for iterative methods such as Newton's iteration are then performed automatically. As all finite elements, including the divergence-conforming ones, are predefined in this software package. The variational equation \eqref{eq:vareq_lag}, or the minimization problem corresponding to \eqref{eq:minimization_deltat} are entered symbolically, using the definitions \eqref{eq:defpsi} for the stored energy and \eqref{eq:defphieps} for the regularized dissipation function $\hat \phi_\varepsilon$.

\section{Numerical results}

\subsection{Non-proportional loading}

As a first numerical example, we consider the repolarization of a piezoelectric cube. This benchmark problem was originally proposed by \cite{HuberFleck:2001}, who provide measurement results. In their experiment, a cube was cut from pre-polarized PZT-5H, such that the axes of the cube are at an angle $\theta$ compared to the original polarization direction. Then, the cube is electroded and re-polarized by applying an electric field aligned with the cube's $z$ axis. The change in dielectric displacement at the electrodes was measured with growing prescribed electric field. 

This example has been widely used as a benchmark example for non-proportional loading e.g. by \cite{McMeekingLandis:2002,Landis:2002, ArockiarajanSivakumarSansour:2009,StarkNeumeisterBalke:2016II}. In our computations, we assume the cube to be of size $\SI{5}{\milli\meter}$. We use the material parameters for PZT-5H as provided by \cite{SemenovLiskowskyBalke:2010}, see Table~\ref{tab:materialPZT5H}. The material is pre-polarized at angle $\theta \in [0^\circ, 180^\circ]$. The regularization parameter for the dissipation function was chosen as $\varepsilon=10^{-4}P_0$. The free energy function was not regularized. 

Concerning this initial state, \cite{StarkNeumeisterBalke:2016II} point out that it is not consistent if $\theta \neq 0^\circ, 180^\circ$, as then surface charges of size $\rho_{P} = \Piv\cdot\nv$ occur at the non-electroded boundaries. In this work, we assume one of their remedy proposals, where these charges are modeled as non-moveable surface charges. This way, on the non-electroded boundaries we have constant surface charges of $\Dv \cdot \nv = \rho_P$. 

\begin{figure}
	\begin{center}
		\includegraphics[width=0.9\columnwidth]{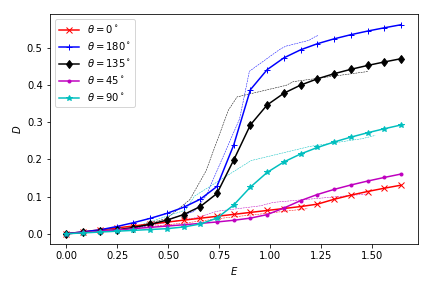}
	\end{center}
\caption{Repoling of piezoelectric cube: The change of dielectric displacement is measured with growing applied electric field for different angles $\theta$. Comparison to the original measurement data taken from \cite{HuberFleck:2001}.}
	\label{fig:repoling_theta}
\end{figure}

In a first computation, we perform repolarization in 20 load steps, where a total electric field of $2 E_0$ is applied. We discretize the cube using an unstructured tetrahedral mesh consisting of 745 elements. We choose the finite element order $k=1$, which means second order displacement elements and first order dielectric/remanent polarization elements. In total, we arrive at 18733 degrees of freedom.
In Figure~\ref{fig:repoling_theta} the change in dielectric displacements as measured from the finite element model is provided, together with the original data as taken from \cite{HuberFleck:2001}. A good qualitative resemblance is observed. 

%Ec = 820e3, Psat = 0.24, eps = 2.77e-8, Emod = 61e9, nu = 0.31, 
%dij = [-2.74e-10, 5.93e-10, 7.41e-10], m=1.4, h0=714000, 
%mc = 1.4, mt = 1.4, h0T = 620e6, T0 = 24.6e6,
%SSat = 9.3e-3, scale=scale, scale_barPi = scale_barPi, is_constepsT=constepsT )

\begin{table}
	\begin{center}
		\begin{tabular}{ll}
			Young's modulus $Y_E$ & $\SI{61e9}{\newton\per\meter\squared}$\\
			Poisson ratio $\nu$ & 0.31\\
			el. permittivity $\epsilon^{S}$ & $\SI{2.77e-8}{\farad\per\meter}$\\
			piezoelectric coupling $d_{31}$ & $\SI{2.74e-10}{\meter\per\volt}$\\
			piezoelectric coupling $d_{33}$ & $\SI{5.93e-10}{\meter\per\volt}$\\
			coercive electric field $E_0$ & $\SI{820e3}{\volt\per\meter}$\\
			saturation polarization $P_0$ & $\SI{0.24}{\coulomb\per\meter\squared}$\\
			saturation strain $S_0$ & $\SI{9.3e-3}{}$\\
			shape parameter $m$ & 1.4\\
			hardening parameter $h_0$ & $\SI{714e3}{\meter\per\farad}$
		\end{tabular}
	\end{center}
\caption{Material parameters for PZT-5H used for repolarization, \cite{SemenovLiskowskyBalke:2010}.} \label{tab:materialPZT5H}
\end{table}

\begin{figure}
\begin{center}
	\includegraphics[width=0.9\columnwidth]{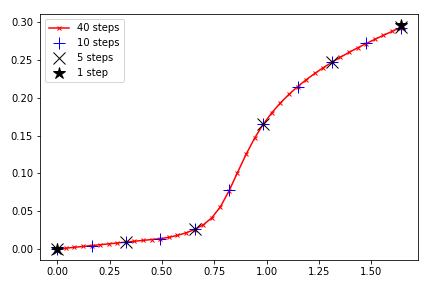}
\end{center}
\caption{Repoling of piezoelectric cube: Comparison of different loadstep sizes for an initial polarization angle of $90^\circ$.}
\label{fig:repoling_stepsize}
\end{figure}

In a second computation, we analyze how the loadstep size affects the results. We consider now the case of $\theta = 90^\circ$. We observe that even large load increments lead to accurate results, see Figure~\ref{fig:repoling_stepsize}. In Table~\ref{tab:relerrors_repoling}, the relative error of the dielectric displacement difference measured at an applied electric field of $2 E_0$ is provided for different numbers of load steps. To this end, the result from applying 40 loadsteps was used as a reference. We see that using only one load step in this highly nonlinear problem led to a relative error below one percent. In this case, a total number of 35 Newton iterations with line search were needed.

\begin{table}
	\begin{center}
		\begin{tabular}{rr}\hline
			number of loadsteps & rel. error\\\hline
			20 & 0.0001241\\
			10 & 0.0005419\\
			5 &  0.0018539\\
			1 &  0.0099080\\\hline
		\end{tabular}
	\end{center}
\caption{Relative error measured for the dielectric displacement difference at fully applied electric field $2 E_0$ for different load step sizes.}
\label{tab:relerrors_repoling}
\end{table}

\subsection{Hexahedron with cylindrical hole}

This example is a benchmark example taken from \cite{SemenovLiskowskyBalke:2010}. A hexahedron of size $\SI{20}{\milli\meter}\times\SI{20}{\milli\meter}\times\SI{6}{\milli\meter}$ made of PZT-5H sports a cylindrical hole of diameter $\SI{4}{\milli\meter}$. Two opposite faces of the hexahedron are electroded, and an electric potential is applied. Material properties are the same as in the previous example and collected in Table~\ref{tab:materialPZT5H}.  The regularization parameter for the dissipation function was chosen as $\varepsilon=10^{-4}P_0$. The free energy function was not regularized. Symmetry of the problem is used, such that the geometry is reduced to one eighth. In the original reference, the potential was applied in five load steps, and iteration counts were provided. We use a prismatic mesh with one prism over the height, which was also done in the original reference. We provide iteration counts for the first and second order method in Table~\ref{tab:itcount}. These counts are higher than those listed by \cite{SemenovLiskowskyBalke:2010}, however we only needed load sub-steps in the second order method in the last load increment. Remanent polarization and strain at $\Delta V = \SI{15}{\kilo\volt}$ are depicted in Figure~\ref{fig:15kv}.

\begin{table}
	\begin{tabular}{rrr}\hline
		&order $k=1$ & order $k=2$\\\hline
		$0\to\SI{5}{\kilo\volt}$ & 8& 10 \\
		$5\to\SI{8}{\kilo\volt}$& 19 & 31\\
		$8\to\SI{9}{\kilo\volt}$& 33 & 34\\
		$9\to\SI{10}{\kilo\volt}$& 21 & 21\\
		$10\to\SI{15}{\kilo\volt}$ & 48 & 13+11+8+7\\\hline
		\end{tabular}
	\caption{Hexahedron with cylindrical hole: Iteration counts in Newton's method with linesearch for the five load steps.}\label{tab:itcount}
\end{table}  

\begin{figure}
	\includegraphics[width=\columnwidth]{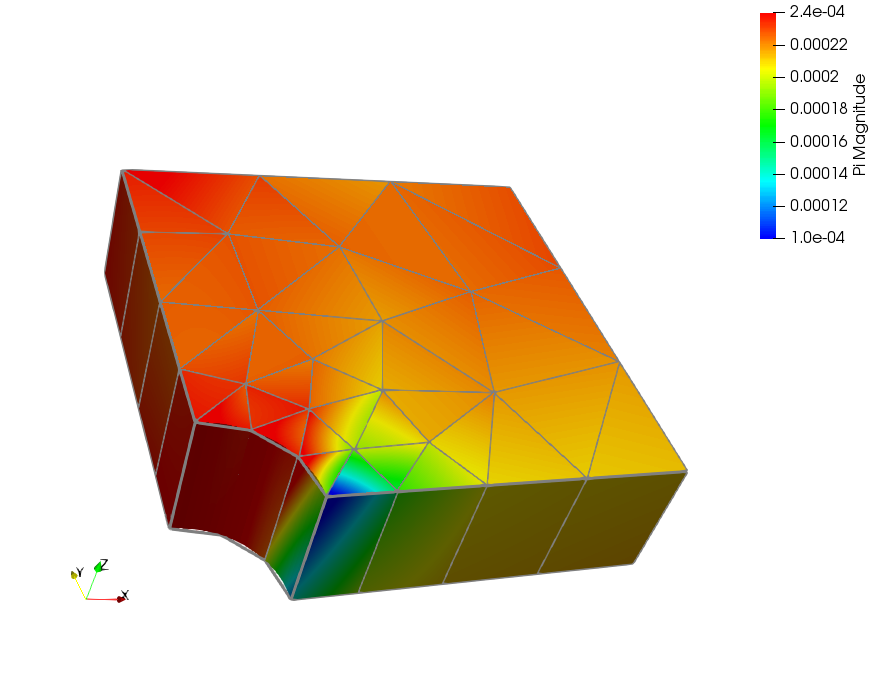}
	\includegraphics[width=\columnwidth]{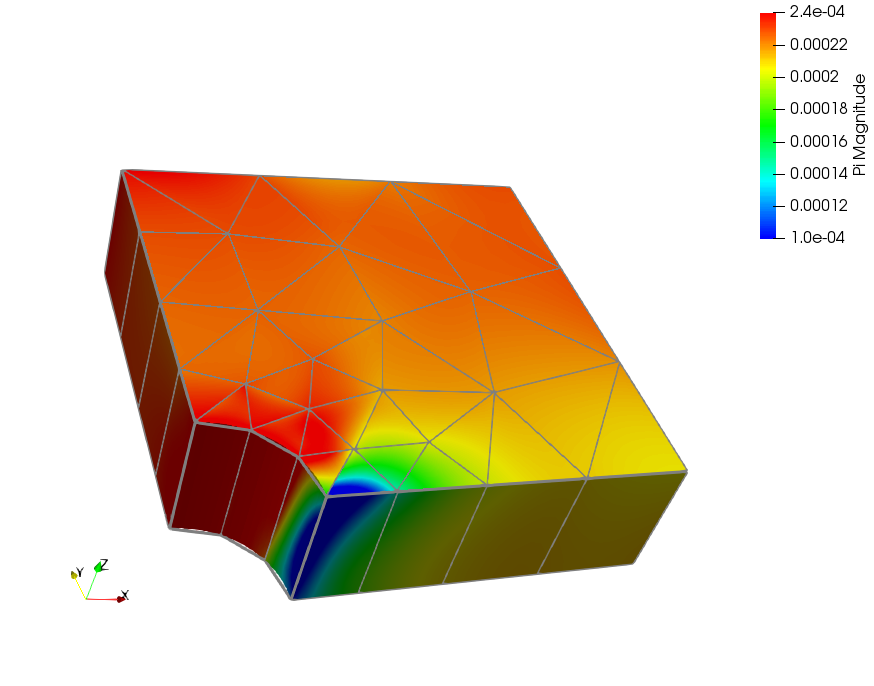}
	\caption{Hexahedron with cylindrical hole: Remanent polarization $|\Piv|$ and strain component $S_{zz}$ at $\Delta V = \SI{15}{\kilo\volt}$.} \label{fig:15kv}
\end{figure}
\begin{figure}
	\includegraphics[width=\columnwidth]{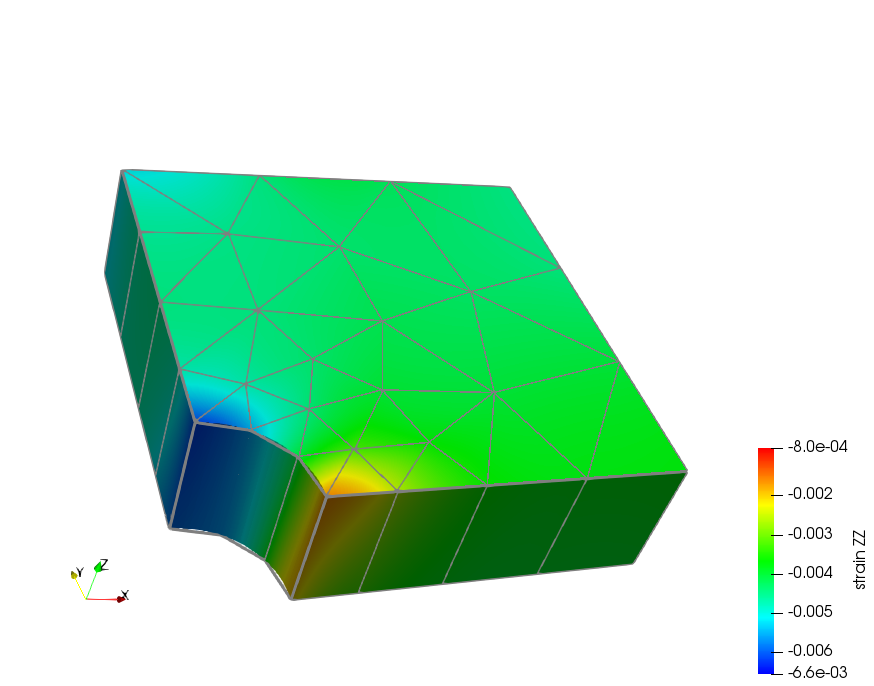}
	\includegraphics[width=\columnwidth]{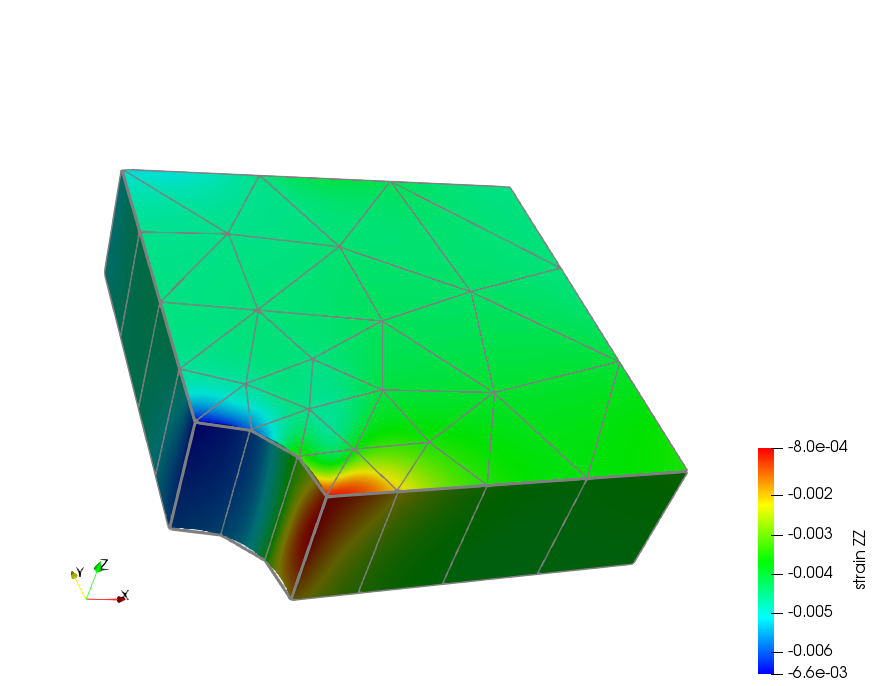}
	\caption{Hexahedron with cylindrical hole: Remanent polarization $|\Piv|$ and strain component $S_{zz}$ at $\Delta V = \SI{15}{\kilo\volt}$.} \label{fig:15kv}
\end{figure}

\subsection{Poling of a $d_{33}$ macro fiber composite}

As an advanced numerical example, we consider a $d_{33}$ macro fibre composite ($d_{33}$ MFC). Such a $d_{33}$ MFC is commercially available e.g. from Smart Materials Corp., Sarasota, FL. This MFC consists of thin (microscopic) fibres of PZT which are aligned in an epoxy matrix. Interfingering copper electrodes are aligned at the top and bottom at an angle of $90^\circ$. The computation of macroscopic material constants resembling the mircostructure has been addressed in the literature in different contributions. Analytic mixing rules for a simplified setup were provided by \cite{DeraemaekerEtal:2009}. In \cite{DeraemaekerNasser:2010}, the authors address the fact that the geometric setup of the $d_{33}$ MFC leads to non-uniform polarization directions in three dimensions. However, they did not simulate the polarization process, but assumed the remanent polarization to be of uniform absolute value $P_0$, but in direction of the electric field. \cite{ShindoEtal:2011} presented a first simulation of the polarization process, but using a different model from the one proposed in this work, and only in two space dimensions. Our contribution is now a complete simulation of the poling process, where a poling electric field of size $4 E_0$ is applied to the interdigitated electrodes. Afterwards, using this remanent polarization field, a small voltage is applied. Strain and dielectric displacement computed in this linear problem then lead to macroscopic $d_{33}$, $d_{32}$ and $\epsilon^T_{33}$ constants.

\begin{figure}
	\begin{center}
		\includegraphics[width=0.98\columnwidth]{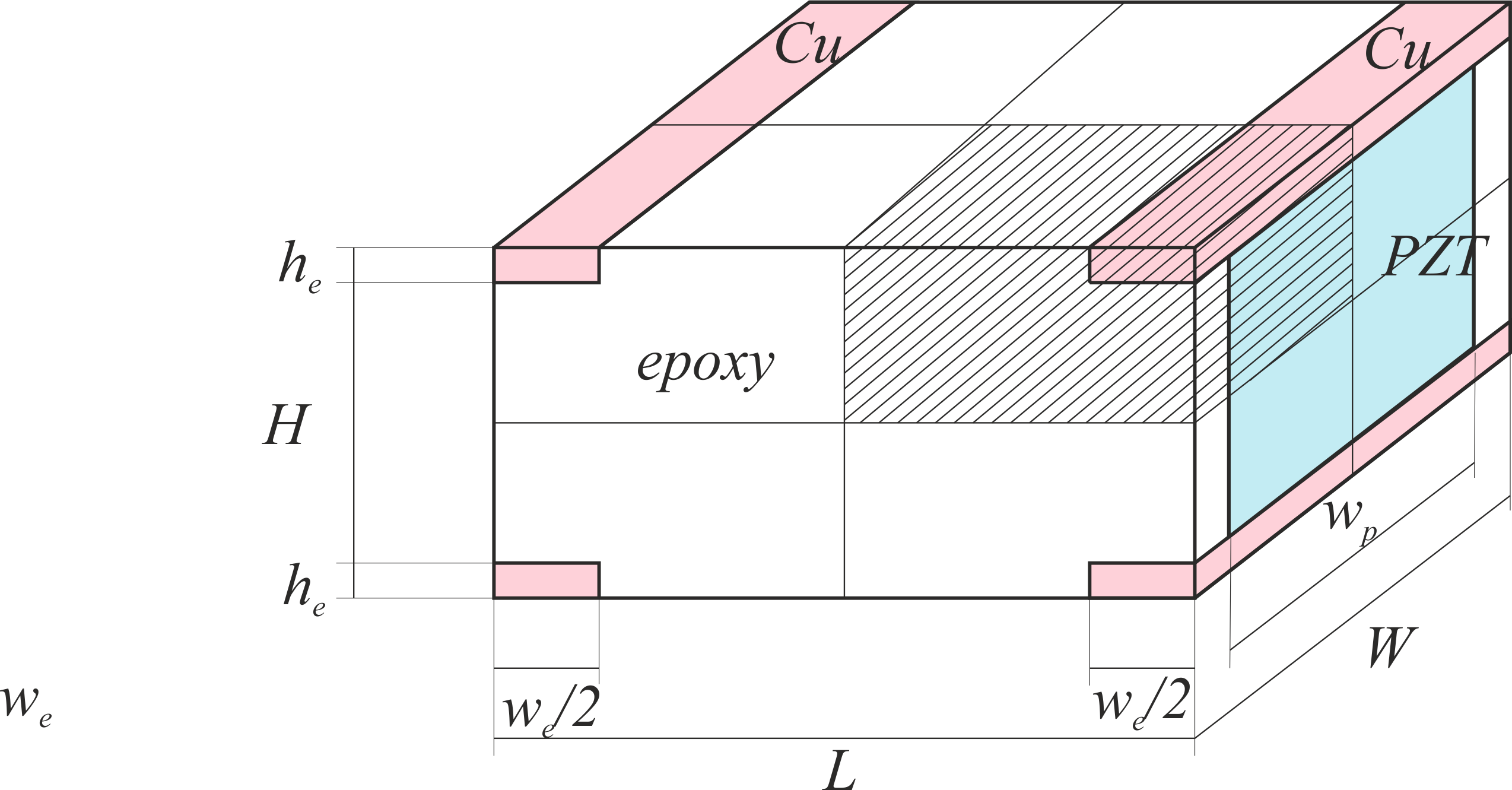}\\
	\end{center}
	\caption{Homogenization of macro fibre composite: schematic sketch of unit cell.}
	\label{fig:MFC_sketch}
\end{figure}

A sketch of the unit cell as used by \cite{DeraemaekerNasser:2010} is provided in Figure~\ref{fig:MFC_sketch}. Conforming with \cite{ShindoEtal:2011}, we use $W = \SI{413.90}{\micro\meter}$ for the width of an epoxy/PZT fibre set. The width of the PZT fibre $w_p$ is then defined by the volume ratio $\rho$ of the fibre. We use $\rho = 0.8$, which implies $w_p = \rho W = \SI{331.12}{\micro\meter}$. Height of fibre and electrode are also taken from \cite{ShindoEtal:2011}, and set to $h_p = \SI{206.40}{\micro\meter}$ and $h_e = \SI{18}{\micro\meter}$. To compare our results to the values provided by \cite{DeraemaekerNasser:2010}, we use their settings $w_e = h_p$ for the width, and $L/w_e = 6$ for the distance of the copper electrodes. This implies $w_e = \SI{206.40}{\micro\meter}$ and $L = \SI{1238.40}{\micro\meter}$.

\begin{figure}
	\begin{center}
		\includegraphics[width=0.98\columnwidth]{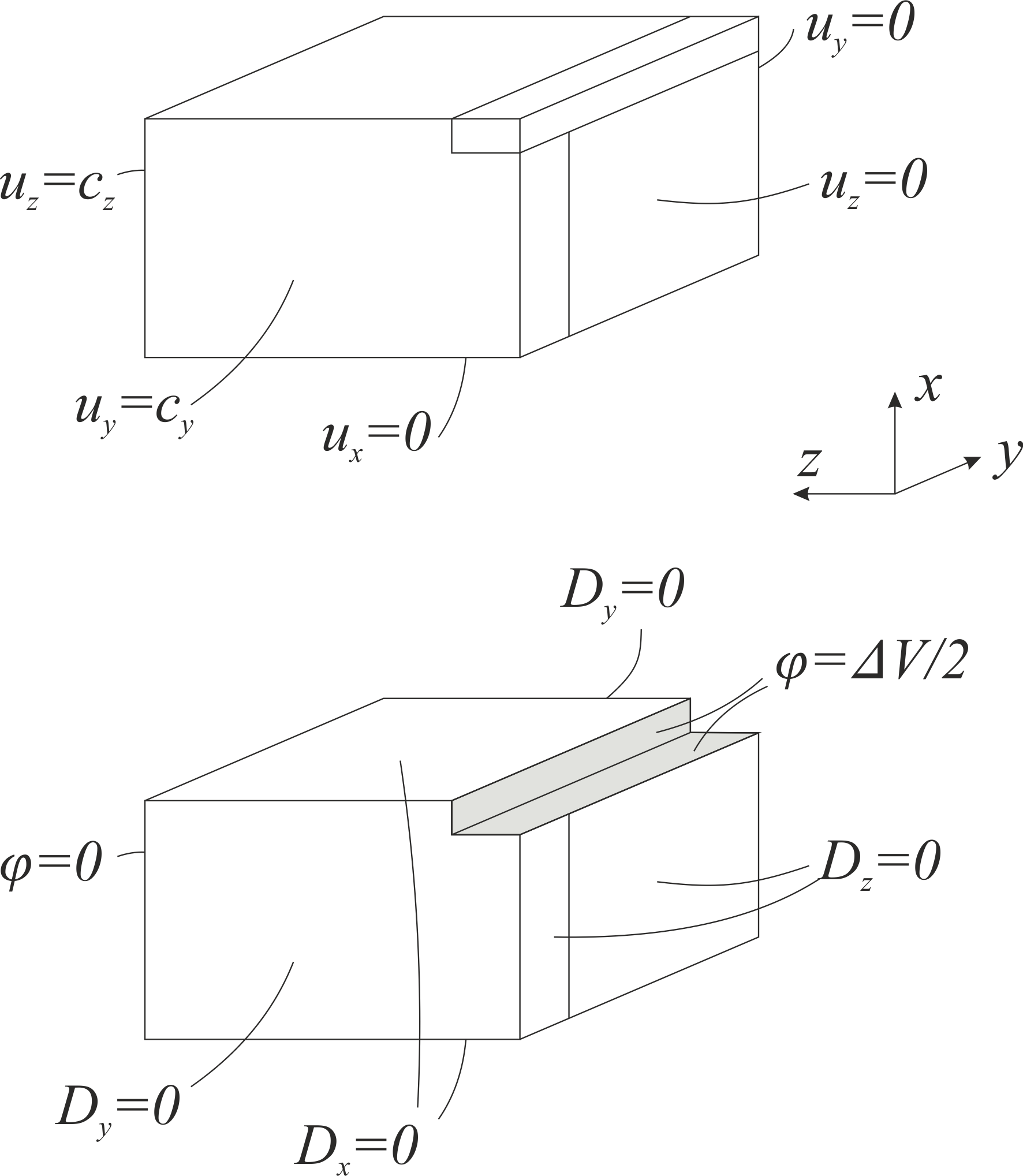}\\
	\end{center}
	\caption{Homogenization of macro fibre composite: symmetry boundary conditions posed on eigth of unit cell, top: displacement boundary conditions, bottom: electric boundary conditions.}
	\label{fig:MFC_bc}
\end{figure}

Due to symmetry of the geometry and the applied load, which is an applied voltage to the electrodes, we may consider an eigth of the geometry only. Boundary conditions are defined in Figure~\ref{fig:MFC_bc}. While the displacement boundary conditions are straightforward for this symmetric problem, we shortly comment on the electric boundary conditions. On the copper electrodes, we assume a constant electric potential of $\varphi = \pm \Delta V$. This is realized in the following way in our mixed finite element formulation: dielectric displacements are not modelled in the electrode volume, but only in the PZT and epoxy parts. On the electrode surfaces, which are now surfaces of the electric domain, the potential boundary condition $\varphi = \pm\frac12 \Delta V$ is applied via a surface integral in the external work statement \eqref{eq:externalwork}. The symmetry plane between electrodes is considered as grounded due to symmetry of the electric potential, which means a (zero) volume integral is added in theory. On all other surfaces, surface charges $\Dv\cdot\nv$ are set to zero.

The choice of material parameters was conducted in the following way: All constants that could be re-used from the linear model by \cite{DeraemaekerEtal:2009,DeraemaekerNasser:2010} were used. This includes the linear properties of epoxy ($Y_E = \SI{2.9}{\giga\pascal}$, $\nu = 0.3$, $\epsilon=4.25 \epsilon_0$) and copper ($Y_E = \SI{117}{\giga\pascal}$, $\nu = 0.31$). In this example, we assume that the stiffness at constant electric field $\cten^{E}$ is constant and isotropic, and defined by Young's modulus $Y_E$ and Poisson ratio $\nu$. The permittivity at constant strain $\epsten^S$ is assumed isotropic, but is not constant. In linear theory, it is smaller than the permittivity at constant stress $\epsten^T$. In this contribution, it is assumed to be a convex combination of $\epsilon^T = 1900 \epsilon_0$ for zero polarization, and $900 \epsilon_0$ for fully polarized material. Further, the additive part of the free energy function $\hat \psi^{i}$ is chosen such that its derivative needed in implementations reads
\begin{align}
\hat \psi^{i}(\Piv) &= \frac{h_0 P_0}{2 |\Piv|} \left(\log\left(\frac{P_0+|\Piv|}{P_0-|\Piv|}\right) \right) \Piv.
\end{align}
The saturation polarization $P_0$ and hardening parameter $h_0$ were chosen such that the hysteresis measurements depicted in \cite[Figure~7(d)]{Hooker:1998} could be reproduced. The remanent saturation strain $S_0$ was chosen in accordance with the measurements provided by \cite{FanStollLynch:1999}.   The regularization parameter for the dissipation function was chosen as $\varepsilon=10^{-4}P_0$, and to avoid convergence problems at the singularity at the electrode also $\hat \psi^{i}$ was regularized using the same parameter.
 All material parameters are also  provided in Table~\ref{tab:materialPZT5A}. 

\begin{table}
	\begin{center}
		\begin{tabular}{ll}
			Young's modulus $Y_E$ & $\SI{48e9}{\newton\per\meter\squared}$\\
Poisson ratio $\nu$ & 0.41\\
el. permittivity $\epsilon^{T}$ & $1900 \epsilon_0$\\
piezoelectric coupling $d_{31}$ & $\SI{-1.85e-10}{\meter\per\volt}$\\
piezoelectric coupling $d_{33}$ & $\SI{4.40e-10}{\meter\per\volt}$\\
coercive electric field $E_0$ & $\SI{1150e3}{\volt\per\meter}$\\
saturation polarization $P_0$ & $\SI{0.255}{\coulomb\per\meter\squared}$\\
saturation strain $S_0$ & $\SI{5.5e-3}{}$\\
hardening parameter $h_0$ & $1/(40\,000\epsilon_0)$
		\end{tabular}
	\end{center}
	\caption{Homogenization of macro fibre composite: material parameters for fibres made from PZT-5A.}
	\label{tab:materialPZT5A}
\end{table}


\begin{figure}
	\begin{center}
		\includegraphics[width=0.98\columnwidth]{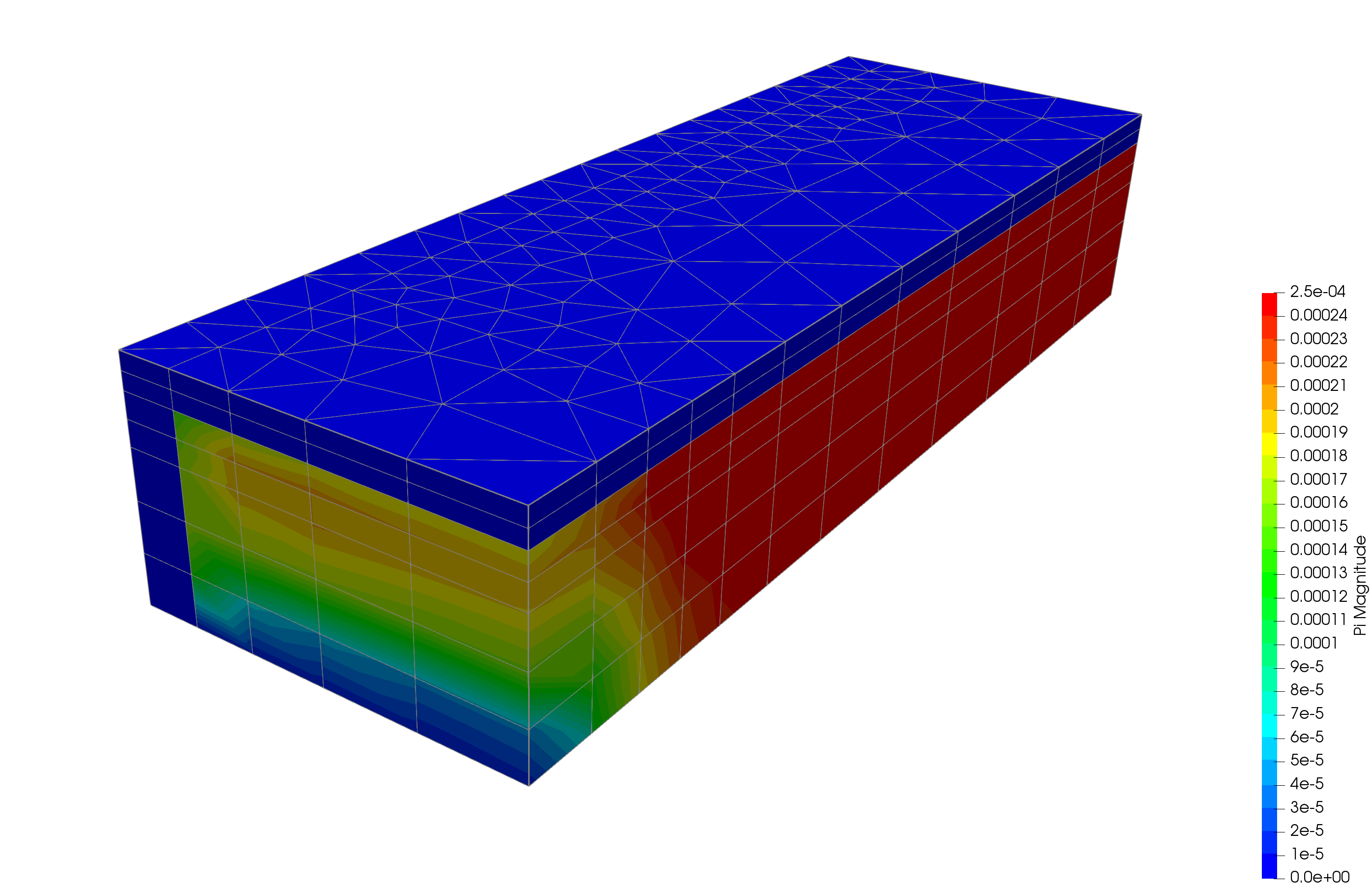}\\
	\end{center}
	\caption{Homogenization of macro fibre composite: remanent polarization at full poling electric field $3E_0$. 
		%Top: fine mesh, finite element order $k=1$, bottom: coarse mesh, finite element order $k=2$.
	}
	\label{fig:MFC_PI}
\end{figure}

\begin{figure}
	\begin{center}
		\includegraphics[width=0.98\columnwidth]{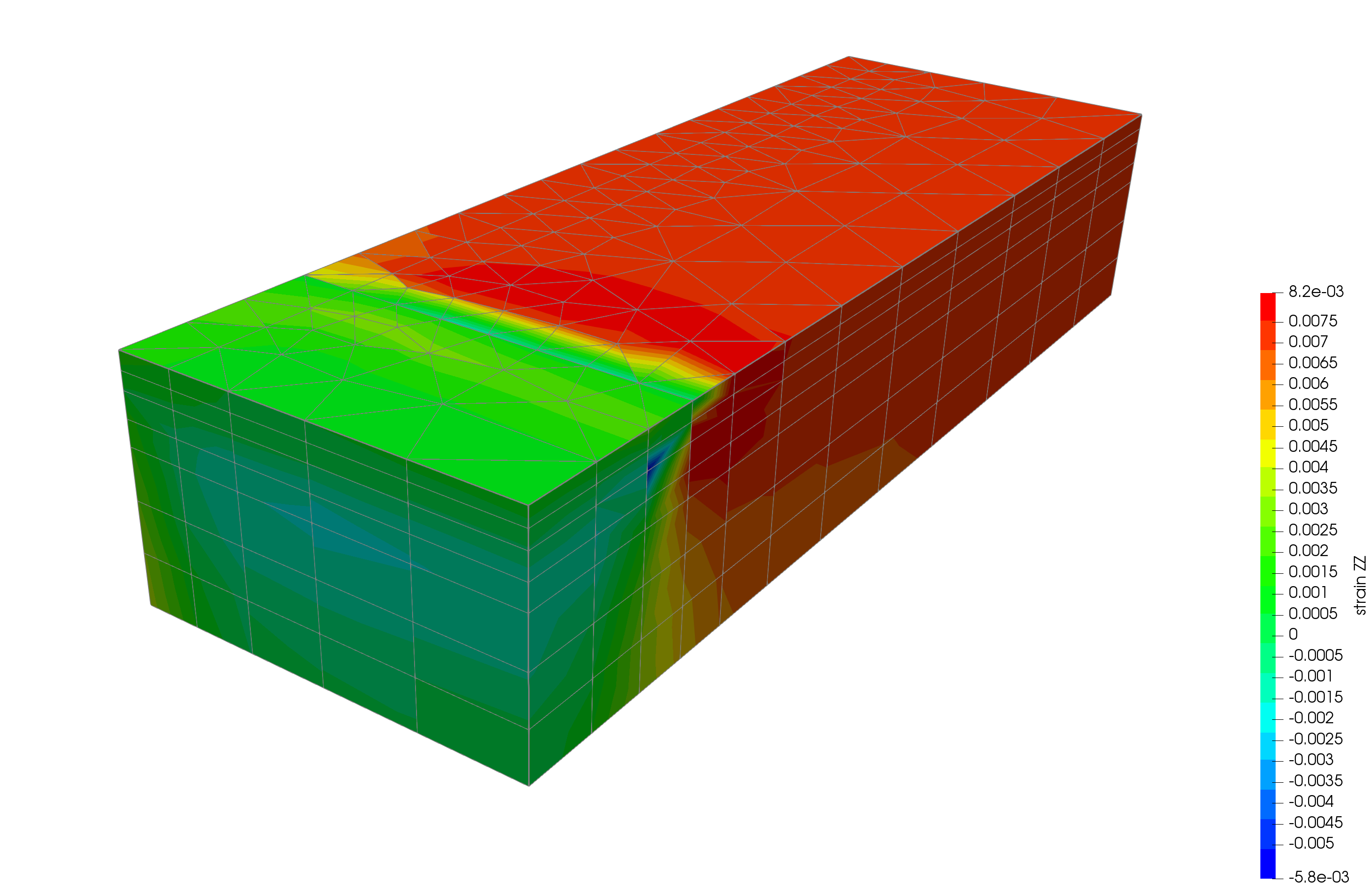}\\
	\end{center}
	\caption{Homogenization of macro fibre composite: longitudinal strain $S_{zz}$ at full poling electric field $3E_0$. 
		%Top: fine mesh, finite element order $k=1$, bottom: coarse mesh, finite element order $k=2$.
	}
	\label{fig:MFC_strainzz}
\end{figure}

\begin{figure}
	\begin{center}
		\includegraphics[width=0.98\columnwidth]{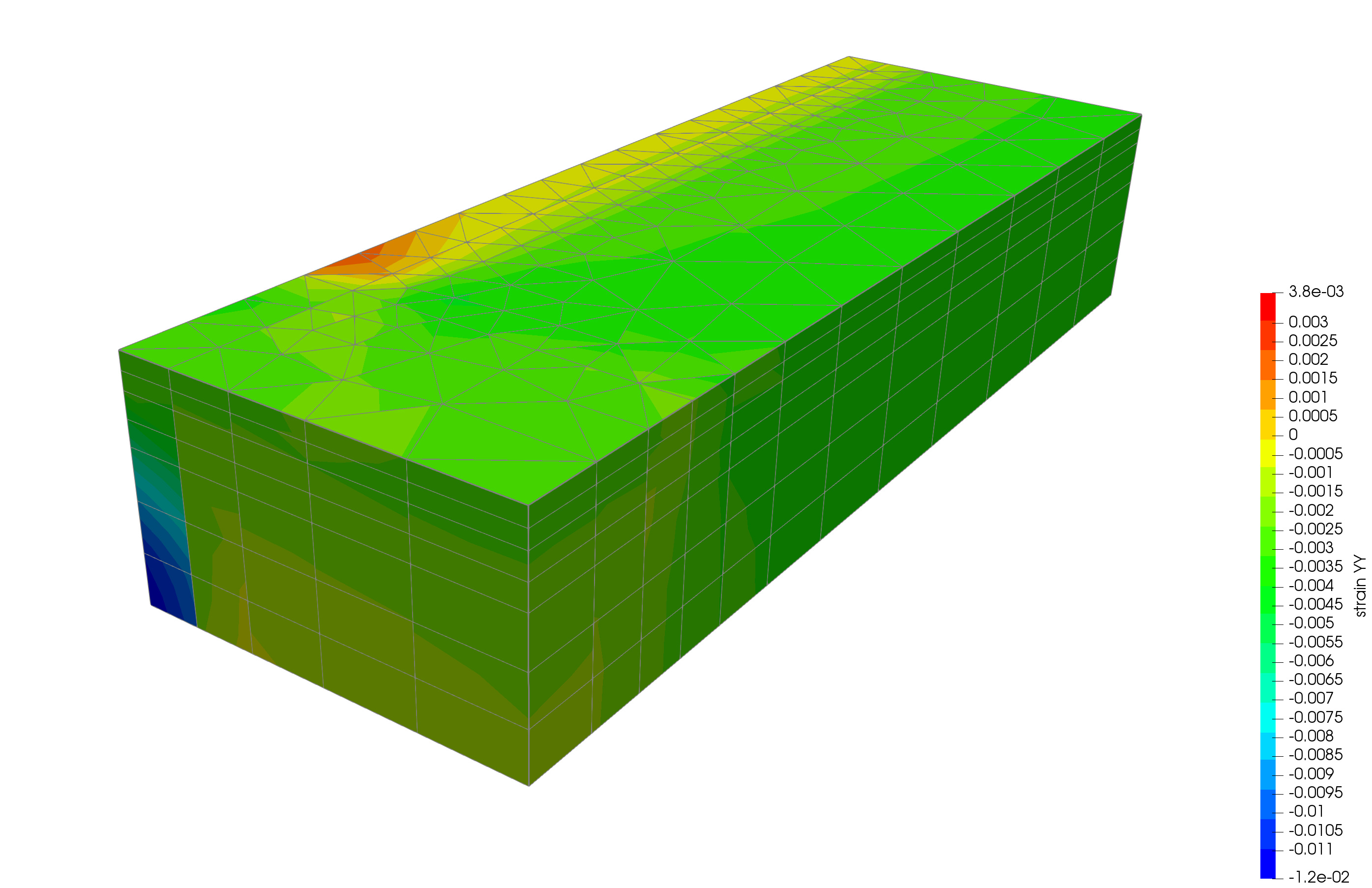}\\
	\end{center}
	\caption{Homogenization of macro fibre composite: transverse strain $S_{yy}$ at full poling electric field $3E_0$. 
		%Top: fine mesh, finite element order $k=1$, bottom: coarse mesh, finite element order $k=2$.
	}
	\label{fig:MFC_strainyy}
\end{figure}

In Figure~\ref{fig:MFC_PI} the absolute value at full poling electric field $3E_0$ is depicted. Figures~\ref{fig:MFC_strainyy} and \ref{fig:MFC_strainzz} show the strain distributions of $S_{yy}$ and $S_{zz}$, respectively. After the poling field was removed, a linear computation was performed to compute the homogenized material constants. To this end, an electric potential was applied to the electrodes, and the macroscopic material constants $d_{32}$, $d_{33}$ and $\epsilon^{T}_{33}$ were computed by the following formulae provided by \cite{DeraemaekerNasser:2010},
\begin{align}
d_{32} &= \frac{\bar S_{yy}}{\bar E_z}, & d_{33} &= \frac{\bar S_{zz}}{\bar E_z}, & \epsilon^{T}_{33} &= \frac{\bar D_z}{\bar E_z}.
\end{align}
The average values of $\bar S_{yy}$ and $\bar S_{zz}$ were derived from the constant displacements arising at the boundaries, compare Figure~\ref{fig:MFC_bc}, $\bar S_{yy} = 2c_y / W$ and $\bar S_{zz} = 2c_z / L$. The average value of the electric field was chosen according to \cite{DeraemaekerNasser:2010} as $\bar E_z = \Delta V/L$. The average value of the dielectric displacement was chosen as the total charge at one electrode divided by the cross section area $W\,H$. Due to Gauss' law and the insulation boundary conditions on all outer surfaces, this value is equal to the average of $\Dv \cdot \nv$ over the symmetric boundary at $L/2$. The following values were found in our computations,
\begin{align}
d_{32} & =  \SI{-1.1554e-10}{\meter\per\volt},\\
d_{33} &=  \SI{3.1598e-10}{\meter\per\volt},\\
\epsilon^T_{33} &=  \SI{1.7736e-08}{\ampere\second\per\volt\per\meter} = 2003 \epsilon_0 .
\end{align}
We note that the piezoelectric constants are lower than those estimated by \cite{DeraemaekerNasser:2010}. There are several facts possibly contributing to this outcome: first, the absolute value of the polarization always stays below saturation polarization. It decreases slightly more as the poling electric field is reduced in the used material law, as the remanent polarization at zero voltage is below the saturation polarization also in one-dimensional settings. Moreover, the polarization is lower close to the fibre boundaries, where the interface to the epoxy matrix is located. And last, in the zones below the electrodes, the electric field and polarization are not oriented in $z$ direction, but rather in vertical $x$ direction. This last issue has also been discussed in the finite element simulation by \cite{DeraemaekerNasser:2010}.

\section{Conclusion}
	In this contribution, we have presented energy-based model of the ferroelectric polarization process. In an energy-based setting, dielectric displacement and strain (or displacement) are the primary independent unknowns, while the remanent polarization vector is added as an internal variable. Electric field and stress are evaluated through post-processing utilizing the material law.  The model is governed by two constitutive functions: the free energy function and the dissipation function. We provide standard choices from the literature for the free energy, and motivate our choice of dissipation function. As it is non-differentiable, we proposed to regularize the problem. Then, a variational equation can be posed, which is subsequently discretized using conforming finite elements for each quantity. One section is devoted to the finite elements, as non-standard continuity is needed for the dielectric displacement to obtain a conforming method. We provide corresponding finite elements, that are chosen such that Gauss' law of zero charges is satisfied exactly. The discretized variational equations are solved for all unknowns at once in a single Newton iteration. We present numerical examples gained in the open source software package Netgen/NGSolve.

\section{Acknowledgements}
	Martin Meindlhumer acknowledges support of Johannes Kepler University Linz, Linz Institute of Technology (LIT).\\
	This work has been supported by the Linz Center of Mechatronics (LCM) in the framework of the Austrian COMET-K2 program.

\end{document}